\documentclass[a4paper,english,10pt, reqno]{amsart}
\usepackage[utf8]{inputenc}
\usepackage[T1]{fontenc}
\usepackage{lmodern}
\usepackage{amssymb}
\usepackage{amsmath}
\usepackage{mathrsfs}
\usepackage{stmaryrd}
\usepackage{enumitem}
\usepackage[olditem,oldenum]{paralist}

\usepackage{graphics}
\usepackage{mdwlist}
\usepackage[all,2cell]{xy} 
\UseAllTwocells
\usepackage{dsfont}
\usepackage{array}
\usepackage{hhline}
\usepackage{multirow}
\usepackage{pdflscape}
\usepackage{rotating}
\usepackage{tabularx}
\usepackage[unicode,bookmarks]{hyperref}
\usepackage[usenames,dvipsnames]{xcolor}
\hypersetup{colorlinks=true,citecolor=NavyBlue,linkcolor=BrickRed,urlcolor=Green} 
\usepackage{calligra}

\usepackage{amsmath, array}
\usepackage{color}

\DeclareMathAlphabet{\mathpzc}{OT1}{pzc}{m}{it}

\DeclareSymbolFont{EulerScripta}{U}{euf}{m}{n}
\SetSymbolFont{EulerScripta}{bold}{U}{euf}{b}{n} 
\DeclareSymbolFontAlphabet\matheufm{EulerScripta}

\DeclareSymbolFont{EulerScriptb}{U}{eur}{m}{n}
\SetSymbolFont{EulerScriptb}{bold}{U}{eur}{b}{n} 
\DeclareSymbolFontAlphabet\matheurm{EulerScriptb}

\DeclareSymbolFont{EulerScriptc}{U}{eus}{m}{n}
\SetSymbolFont{EulerScriptc}{bold}{U}{eus}{b}{n} 
\DeclareSymbolFontAlphabet\matheusm{EulerScriptc}

\usepackage{xspace}

\makeatletter
\@addtoreset{equation}{section}
\renewcommand{\theequation}{\arabic{section}.\arabic{equation}}

\makeatother

 \theoremstyle{plain}
 \newtheorem{theo}{Theorem}[section]
 \newtheorem*{theo*}{Theorem}
  \newtheorem*{prop*}{Proposition}
 \newtheorem{conj}[theo]{Conjecture}
  
 \newtheorem{defi}[theo]{Definition}
 \newtheorem{lemm}[theo]{Lemma}
 \newtheorem{prop}[theo]{Proposition}
 \newtheorem{coro}[theo]{Corollary}

 \theoremstyle{remark}
 \newtheorem{rema}[theo]{Remark}

\newcommand{\appendixref}[1]{\hyperref[#1]{appendix~\ref*{#1}}}
\newcommand{\sectionref}[1]{\hyperref[#1]{section~\ref*{#1}}}
\newcommand{\subsectionref}[1]{\hyperref[#1]{subsection~\ref*{#1}}}
\newcommand{\subsectionsref}[1]{\hyperref[#1]{subsections~\ref*{#1}}}
\newcommand{\theoremref}[1]{\hyperref[#1]{theorem~\ref*{#1}}}
\newcommand{\Theoremref}[1]{\hyperref[#1]{Theorem~\ref*{#1}}}
\newcommand{\lemmaref}[1]{\hyperref[#1]{lemma~\ref*{#1}}}
\newcommand{\remarkref}[1]{\hyperref[#1]{remark~\ref*{#1}}}
\newcommand{\definitionref}[1]{\hyperref[#1]{definition~\ref*{#1}}}
\newcommand{\propositionref}[1]{\hyperref[#1]{proposition~\ref*{#1}}}
\newcommand{\Propositionref}[1]{\hyperref[#1]{Proposition~\ref*{#1}}}
\newcommand{\conjectureref}[1]{\hyperref[#1]{conjecture~\ref*{#1}}}
\newcommand{\corollaryref}[1]{\hyperref[#1]{corollary~\ref*{#1}}}
\newcommand{\Corollaryref}[1]{\hyperref[#1]{Corollary~\ref*{#1}}}
\newcommand{\exampleref}[1]{\hyperref[#1]{example~\ref*{#1}}}
\newcommand{\exerciseref}[1]{\hyperref[#1]{Exercise~\ref*{#1}}}
\renewcommand{\eqref}[1]{\hyperref[#1]{(\ref*{#1})}}
\newcommand{\pararef}[1]{\hyperref[#1]{\S\ref*{#1}}}
\newcommand{\assumptionref}[1]{\hyperref[#1]{Assumption~\ref*{#1}}}
\newcommand{\scholieref}[1]{\hyperref[#1]{Scholie~\ref*{#1}}}
\newcommand{\problemref}[1]{\hyperref[#1]{Problem~\ref*{#1}}}

\usepackage{comment}

\newdir{ >->}{{}*!/-10pt/@{>->}}
\newdir{ >}{{}*!/-10pt/@{>}}
\newdir{ (}{{}*!/-5pt/@^{(}}
\newdir{ )}{{}*!/-5pt/@_{(}}

\newcommand{\pw}{{}^\mathrm{p}\kern-.1pc w}
\newcommand{\cw}{{}^\mathrm{c}\kern-.1pc w}

\newcommand{\cH}{{}^\mathrm{c}\kern-.1pc H}
\newcommand{\pH}{{}^\mathrm{p}\kern-.1pc H}
\newcommand{\ptau}{{}^\mathrm{p}\kern-.1pc \tau}
\newcommand{\ctau}{{}^\mathrm{c}\kern-.1pc \tau}

\newcommand{\sM}{\mathscr M}

\newcommand{\Thom}{\mathrm{Th}}

\newcommand{\bbL}{\mathbf{L}}

\newcommand{\one}{\mathds{1}}
\newcommand{\Hom}{\mathrm{Hom}}
\newcommand{\Spec}{\mathop{\mathsf{Spec}}}

\newcommand{\Osheaf}{\mathcal O}

\newcommand{\scr}{\mathscr}
\newcommand{\Sch}{\mathrm{Sch}}

\newcommand{\SH}{\mathbf{SH}}

\newcommand{\K}{\mathrm{K}}
\newcommand{\Id}{\mathrm{Id}}

\newcommand{\ct}{\mathrm{ct}}

\newcommand{\QUSH}{\mathbf{QUSH}}
\newcommand{\llpar}{(\kern-.15pc(}
\newcommand{\rrpar}{)\kern-.15pc)}

\newcommand{\bbN}{\mathbf{N}}

\newcommand{\bbZ}{\mathbf{Z}}

\newcommand{\bbP}{\mathbf{P}}
\newcommand{\bbV}{\mathbf{V}}
\newcommand{\bbG}{\mathbf{G}}
\newcommand{\ra}{\rightarrow}
\newcommand{\bbA}{\mathbf{A}}

\newcommand{\xra}{\xrightarrow}

\newcommand{\acc}{\breve}

\renewcommand{\llbracket}{[\kern-.08pc[}
\renewcommand{\rrbracket}{]\kern-.08pc]}
\newcommand{\Var}{\mathbf{Var}}

\newcommand{\llangle}{\langle\kern-.12pc\langle}
\newcommand{\rrangle}{\rangle\kern-.12pc\rangle}

\makeatletter
\let\@@seccntformat\@seccntformat
\renewcommand*{\@seccntformat}[1]{%
  \expandafter\ifx\csname @seccntformat@#1\endcsname\relax
    \expandafter\@@seccntformat
  \else
    \expandafter
      \csname @seccntformat@#1\expandafter\endcsname
  \fi
    {#1}%
}
\newcommand*{\@seccntformat@section}[1]{%
  {{\csname the#1\endcsname.}}
}
\newcommand*{\@seccntformat@subsection}[1]{%
  {\bf{\csname the#1\endcsname.}}
}
\newcommand*{\@seccntformat@subsubsection}[1]{%
  {\bf{\csname the#1\endcsname.}}
}

\def\section{\@startsection{section}{1}%
  \z@{.7\linespacing\@plus\linespacing}{.5\linespacing}%
  {\normalfont\scshape\centering}}
\def\subsection{\@startsection{subsection}{2}%
 \z@{.7\linespacing\@plus\linespacing}{.5\linespacing}%
  {\normalfont\bfseries}}
\makeatother

\title[A proof of the integral identity]{A proof of the integral identity via Braden's theorem}
\author{Florian Ivorra}
\address{Univ Rennes, CNRS, IRMAR - UMR 6625, F-35000 Rennes, France}
\email{florian.ivorra@univ-rennes.fr}
\date{Mars 2023}

\begin{document}

\subjclass{14C15, 14F42, 32S30}
\keywords{Nearby motivic sheaves, motivic zeta functions}

\begin{abstract}
The purpose of this paper is to provide a very short proof  of a generalized categorified version, within the motivic stable homotopy category of Morel and Voevodsky, of the integral identity for virtual motives conjectured by Kontsevich and Soibelman. Our proof is an application of an important result in geometric representation theory due to Braden and known as the hyperbolic localization/restriction theorem.  Though originally proved in the context of \'etale sheaves (or sheaves on the associated complex analytic space in the case of complex algebraic varieties) Braden's theorem turns out to hold also in the context of motivic sheaves, at least in the special case of vector bundles with a linear $\bbG_m$-action.
\end{abstract}
\maketitle

\noindent {\bf Acknowledgment.} In an indirect way the present work is a follow-up of \cite{Inventiones} and the joint works with Julien Sebag that led to \cite{Inventiones,AJM}. It benefited greatly from the numerous discussions we had during the writing of those works.

\section{Introduction}

Let $k$ be a field of characteristic zero. In this work we give a short proof  of a generalized categorified version, within the motivic stable homotopy category of Morel and Voevodsky, of the integral identity for virtual motives conjectured by Kontsevich and Soibelman in \cite{KSDT} and proved in a weaker form by  L\^e Quy Thuong  in \cite{MR3299104} and in a stronger form by Nicaise-Payne in \cite{MR3983293}.

Since the work of Behrend \cite{MR2600874} and the subsequent developments, Grothendieck rings of varieties have played a crucial role in the search for motivic refinements of the classical numerical Donaldson-Thomas invariants. In particular, the relation between the Behrend function and vanishing cycles  has given an important role to play in this context to the virtual nearby cycles introduced by Denef-Loeser in \cite{MR1618144} via the use of motivic integration. This interplay between motivic integration and Donaldson-Thomas theory has led as a consequence to conjectural formulas for the computation of virtual nearby cycles in situations pertinent to Donaldson-Thomas theory. The integral identity conjectured by Kontsevich and Soibelman is an example of such a formula which focuses on the computation of virtual nearby cycles of $\bbG_m$-invariant functions on vector bundles with a linear action of $\bbG_m$. 

Though different in nature, both the proof by L\^e Quy Thuong and the proof by Nicaise-Payne work within the framework of motivic integration and make use of the tools it provides. In \cite{MR3299104} Hrushovski-Kazhdan's version of motivic integration is used while in \cite{MR3983293} the integral identity is obtained as an application of a motivic Fubini theorem for the tropicalization map based on Hrushovski-Kazhdan's theory of motivic volumes of semi-algebraic sets.

In \cite{Inventiones} motivic homotopy theory is used instead of motivic integration to prove a generalized categorified version of a formula motivated by Donaldson-Thomas theory and conjectured (proved in special cases) by Behrend-Bryan-Szendr\H oi in \cite{MR3032328} and Davison-Meinhardt in \cite{MR3416109,MR3687097}. The formula focuses on the computation of virtual nearby cycles of weighted $\bbG_m$-equivariant functions on algebraic $k$-varieties endowed with a (multiplicative) action of $\bbA^1$.

Our approach here relies also on motivic stable homotopy and our proof is a very short and immediate application of Braden's hyperbolic localization theorem (in the special case of a vector bundle with a linear $\bbG_m$-action) which though originally proved in \cite{MR1996415} in the context of \'etale sheaves (or sheaves on the associated complex analytic space in the case of complex algebraic varieties) turns out to hold, with the same proof, also in the context of motivic sheaves. 

In comparison with the use of motivic integration, motivic homotopy provides a functorial isomorphism of quasi-unipotent motives, a framework which is much more suited for the study of the monodromy operator than the classical Grothendieck rings of varieties. Indeed, contrary to the monodromic Grothendieck ring of varieties that keeps only track of the semi-simple part of the monodromy, the category of quasi-unipotent motives sees both the semisimple part and the nilpotent part of the monodromy. It also relates directly the integral identity with Braden's theorem which plays a very important role in geometric representation theory. This link potentially extends the range of geometric situations in which an analog of the integral identity might be expected. Indeed Braden's theorem holds in the context of constructible sheaves, $\mathscr D$-modules or \'etale sheaves not only for vector bundles with a linear action of $\bbG_m$ but for any $k$-variety with a $\bbG_m$-action \cite{MR3200429,MR1996415,MR3912059}. It would be interesting to have, if possible, such an extension of Braden's theorem in motivic homotopy theory (the corresponding integral identity would then directly follows from the present work).

\par\smallskip{\itshape{Content of this work}}\par\smallskip

Let us briefly described the organization of this paper. In \sectionref{sec:recollection} we provide a brief overview of the motivic nearby sheaf functor introduced by Ayoub in \cite{AyoubII} and its monodromic counterpart from \cite{AJM} which relies on the category of quasi-unipotent motives defined by Ayoub in \cite{AyoubRigid}. We recall the results from \cite{AJM} that we need and give a brief overview of some consequences of the weak tannakian formalism of Ayoub applied to quasi-unipotent motives \cite{AyoubCrelleI,AyoubCrelleII}.

In \sectionref{sec:mainresults} we state and prove our main result (\theoremref{theo:maintheo}) and give some consequences (such as a functorial generalization of Nicaise-Payne's theorem \cite[Theorem 4.2.1]{MR3983293}). The proof is a direct and short application of Braden's hyperbolic localization theorem in the special case of vector bundle with a linear $\bbG_m$-action. For the sake of completeness, later reference or to convince readers that may have doubt that Braden's proof in its original context still apply to motivic homotopy theory, we have recalled the arguments in an appendix (\theoremref{theo:BradenIso}). 

In \sectionref{sec:virtual} we apply our main result to virtual motives and prove a generalized form of the identity formula in the Grothendieck ring of quasi-unipotent motives (\corollaryref{coro:IIquasiunipotent}).

\par\smallskip{\itshape{Notation and convention}}\par\smallskip

We let $k$ be a field of characteristic zero. A $k$-variety is a separated $k$-scheme of finite type. For simplicity, we write $\bbA^1,\bbG_m$ for $\bbA^1_k,\bbG_{m,k}$ respectively and in the same spirit, we write $-\times-$ for the fiber product $-\times_k-$ over $\Spec(k)$. We denote by $j:\bbG_m\hookrightarrow \bbA^1$ the open immersion by $i:$ the closed immersion. We let $p:\bbA^1\ra \Spec(k)$ and $q:\bbG_m\ra\Spec(k)$ be the projections. 

If $S$ is a $k$-variety, we set $j_S=j\times\Id_S$, $i_S=i\times \Id_S$, $q_S:q\times\Id_S$ and $p_S=p\times\Id_S$.
Given a morphism of $k$-varieties $f:T\ra S$, we set $\acc{f}=f\times\Id_{\bbG_m}$.

If $S$ is a $k$-variety, we consider the Grothendieck ring $\mathscr M_S:=\K_0(\Var_S)[\bbL^{-1}]$ where $\bbL=[\bbA^1\times S]$ and its monodromic variant $\mathscr M^{\hat{\mu}}_S$ (see e.g. \cite[\S2.3]{MR2219263} or \cite[\S5.1]{AJM}). A morphism of $k$-varieties $f:T\ra S$ induces two morphism of rings
\[f^\star:\mathscr M^{\hat{\mu}}_S\ra \mathscr M^{\hat{\mu}}_T\quad\textrm{and}\quad \int_f:M^{\hat{\mu}}_T\ra \mathscr M^{\hat{\mu}}_S\]
defined (here we forget the action of $\hat{\mu}$ to simplify) by $f^\star[Z\ra S]=[Z\times_ST\ra T]$ if $Z\ra S$ is a $S$-variety and $\int_f([a:Y\ra T])=[f\circ a:Y\ra S]$ if $a:Y\ra T$ is a $T$-variety. If there is  no confusion possible $f^\star$ is simply denoted by $-|_T$ and we set $\int_f=\int_T$ if $S=\Spec(k)$ and $f$ is the structural morphism. 
The motivic stable homotopy category of $S$ is denoted by $\SH(S)$ while $\SH_{\ct}(S)$ stands for its full subcategory of compact objects.

\section{Quasi-unipotent motives and nearby motivic sheaves}\label{sec:recollection}

Let $Y$ be a $k$-variety. Consider the open subscheme $\eta=\bbG_m$ of $\bbA^1$ and the zero section $\sigma=\Spec(k)$. Given a morphism $g:Y\ra\bbA^1$ we will consider the commutative diagram
\begin{equation}\label{diag:deg}
\xymatrix{{Y_\eta}\ar[d]^-{g_\eta}\ar[r]^-{j}\ar@{}[rd]|{\square} & {Y}\ar[d]^-{g} & {Y_\sigma}\ar[d]^-{g_\sigma}\ar[l]_-{i}\ar@{}[ld]|{\square}\\
{\eta}\ar[r] & {\bbA^1} & {\sigma}\ar[l]}
\end{equation}
and the associated nearby motivic sheaf functor in motivic stable homotopy
\[\Psi_g:\SH(Y_\eta)\ra\SH(Y_\sigma)\]
constructed by Ayoub in \cite{AyoubII}. The scheme $Y_\eta$ will be called the generic fiber of the morphism $g$ and the scheme $Y_\sigma$ the special fiber of the morphism $g$. Let $\Delta$ be the category of finite ordinals $\underline{n}=\{0<1<\cdots<n\}$ for $n\in\bbN$ and $\bbN^\times=\bbN\setminus \{0\}$ be ordered by the opposite of the division relation. Let 
\begin{equation}\label{eq:motivicUC}
(\theta^{\mathscr R},p_{\Delta\times\bbN^\times}):(\mathscr R,\Delta\times\bbN^\times)\ra\bbG_{m}
\end{equation}
be the morphism of diagrams of $k$-schemes introduced by Ayoub in \cite[D\'efinition 3.5.3]{AyoubII} and consider the motive $\mathscr U=(p_{\Delta\times\bbN^\times})_\sharp(\theta^{\mathscr R})_*\one_{(\mathscr R,\Delta\times\bbN^\times)}$ in $\SH(\eta)$. 
Then, the nearby motivic sheaf functor is given by 
\[\Psi_g=\chi_g\circ (p_{\Delta\times\bbN^\times})_\sharp(\theta^{\mathscr R}_g)_*(\theta^{\mathscr R}_g)^*(p_{\Delta\times\bbN^\times})^*\]
where $\chi_g=i^*j_*$ and $(\theta_g^{\mathscr R},p_{\Delta\times\bbN^\times}):(\mathscr R_g,\Delta\times\bbN^\times)\ra Y_\eta$
is the morphism of diagrams of $k$-schemes obtained from \eqref{eq:motivicUC} by base change along the morphism $g_\eta:Y_\eta\ra\bbG_m$.

\begin{prop}\label{prop:PsiU}
Let $A\in\SH(Y_\eta)$. Then, the canonical morphism
\begin{equation}\label{eq:PsiU}
\Psi_g(A)\ra\chi_g(A\otimes g_\eta^*\mathscr U)
\end{equation}
is an isomorphism.
\end{prop}

\begin{proof}
Recall that the triangulated category $\SH(Y_\eta)$ is generated by the Tate twists of motives of the form $h_{\eta*}\one_{Z_\eta}$ where $h:Z\ra Y$ is a projective morphism and $Z$ is smooth over $k$ (see the first step of the proof of \cite[theorem 4.1.1]{AJM}). Hence, we may assume $A=h_{\eta*}\one_{Z_\eta}$. In that case we have a commutative diagram
\[\xymatrix{{\Psi_g(h_{\eta*}\one_{Z_\eta})}\ar[r]\ar[dd] & {\chi_g(h_{\eta*}\one_{Z_\eta}\otimes g^*_\eta\mathscr U)}\ar[d]\\
{} & {\chi_g(h_{\eta*}(\one_{Z_\eta}\otimes h^*_\eta g^*_\eta\mathscr U))}\ar[d]\\
{h_{\sigma*}\Psi_{g\circ h}(\one_{Z_\eta})}\ar[r] &{\chi_{g\circ h}(\one_{Z_\eta}\otimes (g\circ h)_\eta^*\mathscr U)}}\]
in which the vertical morphisms are isomorphisms by the property {\bf (SPE2)} of \cite[D\'efinition 3.1.1]{AyoubII} and the projection formula (\cite[Th\'eor\`eme 2.3.40]{AyoubI}).
To conclude it suffices to remark that the lower horizontal morphism is an isomorphism (see \cite[\S4.3]{AIS}).
\end{proof}

Let $S$ be a $k$-variety. We recall here the definition of the category of quasi-unipotent motives on $S$ introduced by Ayoub in \cite{AyoubRigid} (see also \cite{AJM}) and its relation with the ring $\mathscr M^{\hat{\mu}}_S$.

 If $Z$ is a smooth $S$-scheme of finite type, $h$ is an element in $\Osheaf(Y)^\times$,  and $n\in\bbN^\times$ is an integer, we consider the $\bbG_{m,S}$-scheme \[a:Q^{\mathrm{gm}}_{n}(Z,h)=\Spec(\Osheaf_Z[t,t^{-1},V]/\langle V^n-ht\rangle)\ra\bbG_{m,S}=\Spec(\Osheaf_S[t,t^{-1}]).\]
 \par The category $\QUSH(S)$ of quasi-unipotent motives over $S$ is defined to be the smallest full triangulated subcategory of $\SH(\bbG_{m,S})$ stable under all (small) direct sums and containing the objects of the form 
 \begin{equation}\label{eq:defiquasi}
 a_\sharp\one_{Q^{\mathrm{gm}}_n(Z,h)}(r)=a_!a^!\one_{\bbG_{m,S}}(r)
 \end{equation}
 where $r\in\bbZ$ is an integer (see  \cite[D\'efinition 1.3.25]{AyoubRigid} and \cite[\S 3.2]{AJM}). The full subcategory $\QUSH_\ct(S)$ of compact objects in $\QUSH(S)$ coincides with the smallest full triangulated subcategory of $\SH(\bbG_{m,S})$ containing the objects of the form \eqref{eq:defiquasi} and stable under direct summands.

 As shown in \cite[Lemma 2.1]{IS} (see also \cite[Lemma 8.5]{AIS}), there exists a ring morphism 
\[\chi_{S,c}:\sM_S\ra\K_0(\SH_{\ct}(S))\]
which is uniquely defined by the formula $\chi_{S,c}([Y])=g_!\one_Y$ for $g:Y\ra S$ a morphism of $k$-varieties. Note that in particular $\chi_{S,c}(\bbL)=[\one_{S\times\bbG_m}(-1)]$. If $k$ contains all roots of unity, this morphism is refined in \cite[\S5.2]{AJM} to take into account the monodromy. More precisely, a ring morphism
\[\chi_{S,c}^{\hat{\mu}}:\sM^{\hat{\mu}}_S\ra \K_0(\QUSH_\ct(S))\]
is constructed  which makes the square
\[\xymatrix{{\sM^{\hat{\mu}}_S}\ar[d]^-{\textrm{forget}}\ar[r] & {\K_0(\QUSH_\ct(S))}\ar[d]^{1^*_S}\\
{\sM_S}\ar[r] & {\K_0(\SH_\ct(S))}}\]
commutative. Here the morphism $1^*_S$ is the morphism induced by the pullback $1^*_S:\QUSH_\ct(S)\ra\SH_{\ct}(S)$ along the unit section of $S\times\bbG_m$.

\begin{lemm}\label{lemm:chic}
Let $f:T\ra S$ be a morphism of $k$-varieties. Then, the squares
\[\xymatrix{{\sM^{\hat{\mu}}_S}\ar[r]^-{\chi^{\hat{\mu}}_{S,c}}\ar[d]^{f^\star} & {\K_0(\QUSH_\ct(S))}\ar[d]^-{\acc{f}^*}\\
{\sM^{\hat{\mu}}_T}\ar[r]^-{\chi^{\hat{\mu}}_{T,c}} & {\K_0(\QUSH_\ct(T)),}}\qquad
\xymatrix{{\sM^{\hat{\mu}}_T}\ar[r]^-{\chi^{\hat{\mu}}_{T,c}}\ar[d]^-{\int_f} & {\K_0(\QUSH_\ct(T))}\ar[d]^-{\acc{f}_!}\\
{\sM^{\hat{\mu}}_S}\ar[r]^-{\chi^{\hat{\mu}}_{S,c}} & {\K_0(\QUSH_\ct(S))}}\]
are commutative.
\end{lemm}
\begin{proof}
The commutativity of the square on the left-hand side follows immediately from the proper base change theorem while the commutativity of square on the right-hand side is just a consequence of functoriality. 
\end{proof}

Let us now recall the monodromic nearby motivic sheaf functor and the main result of \cite{AJM} that relates quasi-unipotent motives and motivic nearby sheaves. Let $g:Y\ra\bbA^1$ be a morphism of $k$-varieties and $g^{\bbG_m}:Y\times\bbG_{m}\ra\bbA^1$ be the morphism 
\[
Y\times_k\bbG_m=\Spec(\Osheaf_Y[T,T^{-1}]) \ra\bbA^1_k=\Spec(k[t]).
\]
defined by $t\mapsto Tg$. The monodromic nearby motivic sheaf functor $\Psi^{\textrm{mon}}_g$ introduced in \cite[\S4.1]{AJM} is defined by
\[\Psi^{\textrm{mon}}_g:=\Psi_{g^{\bbG_m}}q_{Y_\eta}^*.\]
The diagram \eqref{diag:deg} admits the following factorization
\[
\xymatrix{{Y_\eta}\ar[r]\ar@{}[rd]|{\square}\ar[d]^-{1_{Y_\eta}} & {Y}\ar[d]^-{1_Y}\ar@{}[rd]|{\square} & {Y_\sigma}\ar[l]\ar[d]^{1_{Y_\sigma}}\\
{\bbG_m\times Y_\eta}\ar[d]^-{g^{\bbG_m}_{\eta}}\ar[r]\ar@{}[rd]|{\square} & {\bbG_m\times Y}\ar[d]^-{g^{\bbG_m}}\ar@{}[rd]|{\square} & {\bbG_m\times Y_\sigma}\ar[l]\ar[d]^-{g^{\bbG_m}_{\sigma}}\\
{\eta}\ar[r] & {\bbA^1} & {\sigma}\ar[l]} 
\]
and {\bf{(SPE2)}} of \cite[D\'efinition 3.1.1]{AyoubII} provides a canonical morphism, for every object  $A$ in $\SH(Y_\eta)$,
\[1^*_{Y_\sigma}\Psi^{\textrm{mon}}_g(A)\ra \Psi_g(A) \]
which is functorial in $A$. We know by \cite[Theorem 4.1.1]{AJM} that this morphism is an isomorphism and by \cite[Theorem 4.2.1]{AJM} that $\Psi^{\textrm{mon}}_g(A)$ is a quasi-unipotent motive over $Y_\sigma$. As a consequence of the weak tannakian formalism developed by Ayoub in \cite[\S2.2]{AyoubCrelleII}, the above isomorphism turns $\Psi_g(A)$ into a comodule over the Hopf algebra $1^*\mathscr U$ (see \cite[Corollaire 2.19]{AyoubCrelleII} for the understanding of the Betti realization of such an action).  

Assume that $k$ contains all roots of unity.  In \cite[Theorem 5.3.1]{AJM} we prove the equality
\[\chi_{Y_\sigma,c}(\psi_g)=\left[\Psi_g^{\textrm{mon}}(\one_{Y_\eta})\right]\]
in the ring $\K_0(\QUSH(Y\sigma))$. In this equality the element $\psi_g$ of $\mathscr M^{\hat{\mu}}_{Y_\sigma}$ is the virtual nearby cycles introduced by Denef and Loeser via motivic integration. Let us sketch the definition of  $\psi_g$ (we refer to \cite{MR1618144} or \cite[Chap. 7,\S4]{ACLNS} for details). The jet scheme $\mathscr L_n(Y)$ of level $n$ is the $k$-scheme canonically defined by the following universal property:
\[
\Hom_{\Sch_k}(\Spec(A),\scr L_n(Y))\cong\Hom_{\Sch_k}(\Spec(A[[T]]/\langle T^{n+1}\rangle),Y)
\]
for every $k$-algebra $A$. Let $d$ be the dimension of $Y$. By analogy with the Igusa $p$-adic integrals, Denef and Loeser have introduced the formal power series:
\[
Z(g,T)=\sum_{n\geq 1}[\scr X_n(g)]\mathbf L^{-nd}T^n\in \mathscr M^{\hat{\mu}}_{Y_\sigma}[[T]]
\]
where $\mathscr X_n(g)$ denotes the locally closed subset consisting of the elements $\phi$ in $\mathscr L_n(Y)$ with base point in $Y_\sigma$ and such that $g(\phi(T)=T^n\mod T^{n+1}$.  Using motivic integration, Denef and Loeser have shown that this formal power series is a rational function and defined $\psi_g$ as minus the limit at infinity (in a naive sense) of $Z(g,T)$.

\begin{rema}
Note that both the Grothendieck ring $\mathscr M^{\hat{\mu}}_{Y_\sigma}$ and $\psi_g$ keep track of the monodromy only in a very weak way as they see only the semi-simple part of monodromy action. The loss of the nilpotent part of the monodromy is a very important drawback of the motivic integration approach (aside from the obvious one that we are not dealing with functorial construction that directly relate with cohomology but rather with Euler characteristics of said cohomology).  This problem disappears if one is willing to work in motivic stable homotopy theory with $\Psi_g$ as shown by Ayoub (see \cite[\S3.6]{AyoubII} and \cite[\S11]{AyoubEtale}).
\end{rema}

\section{Main results}\label{sec:mainresults}

Let $k$ be a field of characteristic zero and $X$ be a $k$-variety.
Let $\mathscr E$ be a locally free $\Osheaf_X$-module of finite rank and denote by $\bbV(\mathscr E)$ the associated vector bundle on $X$. Suppose that a linear action 
\[a:\bbG_m\times\bbV(\mathscr E)\ra\bbV(\mathscr E)\]
of $\bbG_m$ on $\bbV(\mathscr E)$ is given. This action corresponds to a morphism $\bbG_m\times X\ra\mathbf{GL}(\mathscr E)$ of group schemes over $X$ and we may consider the associated weight decomposition 
\[\mathscr E=\bigoplus_{n\in \bbZ}\mathscr E_n.\]
We let $\mathscr E^+=\bigoplus_{n\geq 0}\mathscr E_n$ and $\mathscr E^-=\bigoplus_{n\leq 0}\mathscr E_n$ be respectively the nonnegative and nonpositive parts. The scheme of fixed points under the action of $\bbG_m$ is the closed subscheme $\bbV(\mathscr E^0)$ of $\bbV(\mathscr E)$. We have a commutative diagram
\begin{equation}\label{eq:diagnot}
\xymatrix{{\bbV(\mathscr E^0)}\ar[r]^-{s^-}\ar[d]_-{s^+}\ar@{}[rd]|{\square} & {\bbV(\mathscr E^-)}\ar[d]^-{e^-}\ar[rdd]^-{\pi^-} &{}\\
{\bbV(\mathscr E^+)}\ar[r]^-{e^+}\ar[rrd]_-{\pi^+} & {\bbV(\mathscr E)}\ar[rd]^-{\pi}&{}\\
{} & {}&{\bbV(\mathscr E^0)}}
\end{equation}
in which $s^-,e^-,s^+,e^+$ are closed immersions and $\pi^-,\pi^+,\pi$ are the projections.

Let $f:\bbV(\mathscr E)\ra\bbA^1$ be a $\bbG_m$-invariant function and denote by $f^0:\bbV(\mathscr E^0)\ra \bbA^1$ its restriction to the scheme of fixed points. 
Given an object $A\in\SH(\bbV(\mathscr E)_\eta)$, we have a canonical morphism
\[\pi^+_{\sigma!}e^{+*}_\sigma\Psi_f(A)\ra \pi^+_{\sigma!}\Psi_{f\circ e^+}(e^{+*}_\eta A)\]
provided by property {\bf{(SPE2)}} in \cite[Definition 3.1.1]{AyoubII} and a canonical morphism
\[\pi^+_{\sigma!}\Psi_{f^0\circ \pi^+}(e^{+*}_\eta A)\ra\Psi_{f^0}(\pi^+_{\eta!}e^{+*}_\eta A)
\]
constructed in \cite[(11) p.10]{AyoubII}. Since $f$ is $\bbG_m$-invariant, we have $f\circ e^+=f^0\circ\pi^+$ and by composition we get a canonical morphism 
\begin{equation}\label{eq:MainTheo}
\pi^+_{\sigma!}e^{+*}_\sigma\Psi_f(A)\ra \Psi_{f^0}(\pi^+_{\eta!}e^{+*}_\eta A).
\end{equation}
We may now state our main result:
\begin{theo}\label{theo:maintheo}
Let $f:\bbV(\mathscr E)\ra \bbA^1$ be a $\bbG_m$-invariant morphism and $A\in\SH(\bbV(\mathscr E)_\eta)$ be a $\bbG_m$-equivariant object (in the sense of \definitionref{defi:equivariant}). Then, the morphism 
\[
\pi^+_{\sigma!}e^{+*}_\sigma\Psi_f(A)\ra \Psi_{f^0}(\pi^+_{\eta!}e^{+*}_\eta A)
\]
constructed in \eqref{eq:MainTheo} is an isomorphism in $\SH(\bbV(\mathscr E^0)_\sigma)$.
\end{theo}

\begin{proof}
Let $B$ an object in $\SH(\bbV(\mathscr E)_\eta)$.
Using property {\bf{(SPE2)}} in \cite[Definition 3.1.1]{AyoubII}, in a similar way as in the construction of the morphism \eqref{eq:MainTheo}, we obtain a canonical morphism 
\begin{equation}\label{eq:ProofMainTheo}
\pi^+_{\sigma!}e^{+*}_\sigma\chi_f(B)\ra \chi_{f^0}(\pi^+_{\eta!}e^{+*}_\eta B).
\end{equation}
Then, the morphism \eqref{eq:MainTheo} fits into the commutative diagram
\[\xymatrix{{\pi^+_{\sigma!}e^{+*}_\sigma\Psi_f(A)}\ar[dd]^-{\eqref{eq:MainTheo}}\ar[r]^-{\eqref{eq:PsiU}} & {\pi^+_{\sigma!}e^{+*}_{\sigma}\chi_f(A\otimes f^*_\eta\mathscr U)}\ar[r]^-{\eqref{eq:ProofMainTheo}} & {\chi_{f^0}(\pi^+_{\eta!}e^{+*}_\eta (A\otimes f^*_\eta\mathscr U))}\ar@{=}[d]\\
{} & {} & {\chi_{f^0}(\pi^+_{\eta!}(e^{+*}_\eta A\otimes \pi^{+*}_\eta f^{0*}_\eta\mathscr U))}\ar[d]\\
{\Psi_{f^0}(\pi^+_{\eta!}e^{+*}_\eta A)}\ar[rr]^-{\eqref{eq:PsiU}} & {} & {\chi_{f^0}(\pi^+_{\eta!}e^{+*}_\eta A\otimes f^{0*}_\eta\mathscr U)}}\]
in which the unlabeled morphism is the isomorphism induced by the projection isomorphism (see \cite[Th\'eor\`eme 2.3.40]{AyoubI}). Note that $A\otimes f_\eta^*\mathscr U$ is $\bbG_m$-equivariant and that the morphisms \eqref{eq:PsiU} are isomorphisms by \propositionref{prop:PsiU}. Hence, to prove \theoremref{theo:maintheo}, it suffices to show that if $B$ in $\SH(\bbV(\mathscr E)_\eta)$ is $\bbG_m$-equivariant, then the morphism \eqref{eq:ProofMainTheo} is an isomorphism.

For this consider the commutative diagram
\[\xymatrix{{\bbV(\mathscr E)_\sigma}\ar[d]^-i & {\bbV(\mathscr E^+)_\sigma}\ar[d]^{i^+}\ar[r]^{\pi^+_\sigma}\ar[l]_-{e^+_\sigma}  &{\bbV(\mathscr E^0)_\sigma}\ar[d]^-{i^0}\\
{\bbV(\mathscr E)} & {\bbV(\mathscr E^+)}\ar[r]^-{\pi^+}\ar[l]_-{e^+}  &{\bbV(\mathscr E^0)}\\
{\bbV(\mathscr E)_\eta}\ar[u]_-j & {\bbV(\mathscr E^+)_\eta}\ar[u]_-{j^+}\ar[r]^-{\pi^+_\eta}\ar[l]_-{e^+_\eta}  &{\bbV(\mathscr E^0)_\eta.}\ar[u]_-{j^0}}\]
The morphism \eqref{eq:ProofMainTheo} is the composition of the morphism
\[\pi^+_{\sigma!}e^{+*}_\sigma i^*j_*B\ra i^{0*}\pi^+_!e^{+*}j_*B\]
which is an isomorphism by the proper base change theorem and the image by the functor $i^{0*}$ of the morphism
\begin{equation}\label{eq:morchi}
\pi^+_!e^{+*}j_*B\ra j^0_*\pi^+_{\eta!}e^{+*}_\eta B.
\end{equation}
Therefore, it is enough to show that this last morphism is an isomorphism. We have a commutative diagram
 \[\xymatrix{{\pi^-_*e^{-!}j_*B}\ar[r]\ar[d] & {\pi^-_*j^-_*e^{-!}_\eta B}\ar@{=}[r] & {j^0_*\pi^-_{\eta*}e^{-!}_\eta B}\ar[d]\\
 {\pi^+_!e^{+*}j_* B}\ar[rr]^-{\eqref{eq:morchi}} & & {j^0_*\pi^+_{\eta!}e^{+*}_\eta B}}\]
 in which the vertical morphism on the left is an isomorphism by Braden's \theoremref{theo:BradenIso} and the top horizontal morphism is an isomorphism by the smooth base change theorem. Hence, to prove \theoremref{theo:maintheo} we are left to show that the morphism $\pi^-_{\eta*}e^{-!}_\eta B\ra \pi^+_{\eta!}e^{+*}_\eta B$ is an isomorphism. The open immersion $j:\bbV(\mathscr E)_\eta\hookrightarrow\bbV(\mathscr E)$ is $\bbG_m$-equivariant and thus $j_!B$ is a $\bbG_m$-equivariant object in $\SH(\bbV(\mathscr E))$. By applying Braden's \theoremref{theo:BradenIso} we know that the morphism
\[\pi^-_*e^{-!}j_!B\ra \pi^+_!e^{+*}j_!B\]
is an isomorphism in $\SH(\bbV(\mathscr E^0))$. By applying $j^{0*}$ and using the smooth and proper base change theorem, we deduce that the morphism
\[\pi^-_{\eta*}e_\eta^{-!}j^*j_!B\ra \pi^+_{\eta!}e^{+*}_\eta j^*j_!B\]
is an isomorphism. To conclude, it is enough to observe that the morphism $j^*j_!B\ra B$ is an isomorphism in $\SH(\bbV(\mathscr E)_\eta)$.
\end{proof}

In particular \theoremref{theo:maintheo} provides an isomorphism 
\begin{equation}\label{eq:morII}
\pi^+_{\sigma!}e^{+*}_\sigma\Psi_f(\one_{\bbV(\mathscr E)_\eta})\ra \Psi_{f^0}(\pi^+_{\eta!}e^{+*}_\eta \one_{\bbV(\mathscr E)_\eta}).
\end{equation}
This isomorphism can be thought as the expression in motivic stable homotopy theory of (a generalized form of) the integral identity of Kontsevich and Soibelman. If $x$ in a rational point inside the special fiber of $f^0$, then the pullback of \eqref{eq:morII} along the morphism $x:\Spec(k)\ra\bbV(\mathscr E^0)$ provides a categorified version of Nicaise-Payne's theorem \cite[Theorem 4.2.1]{MR3983293} (see \corollaryref{coro:IImonodromic} for the monodromic version). To understand how the isomorphism \eqref{eq:morII} relates to the classical formulations of the integral identity for virtual motives that appear in \cite{KSDT,MR3299104,MR3983293} (see \sectionref{sec:virtual}), it is useful to reformulate the right-hand side of \eqref{eq:morII} using inverse Thom equivalences. 

Let $p:V\ra S$ be a vector bundle on a $k$-variety $S$. Recall from \cite[D\'efinition 1.5.1]{AyoubI} that the inverse Thom equivalence associated with $V$ is the endofunctor of $\SH(S)$ defined by $\mathsf{Th}(s,p)=s^!p^*$ where $s$ denotes the zero section of the vector bundle $V$. Note that if we set $\Thom^{-1}(s,p):=\mathsf{Th}^{-1}(s,p)(\one_S)$, 
then for every object $A$ in $\SH(S)$ one has a natural isomorphism (see \cite[p.395]{AyoubI})
\[\Thom^{-1}(s,p)\otimes A\xra{\simeq}\mathsf{Th}^{-1}(s,p)(A).\]

Since the motives 
$\pi^+_{\eta!}e^{+*}_\eta\one_{\bbV(\mathscr E)_\eta}$ and $\Thom^{-1}(s^+_\eta,\pi^+_\eta)$ are canonically isomorphic, \cite[Proposition 3.1.7]{AyoubII} provides a canonical isomorphism
\begin{equation}\label{eq:reformulation}
\Psi_{f^0}(\Thom^{-1}(s^+_\eta,\pi^+_\eta))\xra{\simeq} \Thom^{-1}(s^+,\pi^+)|_{\bbV(\mathscr E^0)_\sigma}\otimes\Psi_{f^0}(\one_{\bbV(\mathscr E^0)_\eta})
\end{equation}
where $\Thom^{-1}(s^+,\pi^+)|_{\bbV(\mathscr E^0)_\sigma} $ is the inverse image of $\Thom^{-1}(s^+,\pi^+)$ along the closed immersion of $\bbV(\mathscr E^0)_\sigma$ into $\bbV(\mathscr E^0)$.
Hence, as a corollary of \theoremref{theo:maintheo}, we get the following result.

\begin{coro}\label{coro:maincoro}
Let $f:\bbV(\mathscr E)\ra \bbA^1$ be a $\bbG_m$-invariant morphism. Then, the morphism 
\[
\pi^+_{\sigma!}e^{+*}_\sigma\Psi_f(\one_{\bbV(\mathscr E)_\eta})\ra  \Thom^{-1}(s^+,\pi^+)|_{\bbV(\mathscr E^0)_\sigma}\otimes\Psi_{f^0}(\one_{\bbV(\mathscr E^0)_\eta})
\]
obtained as the composition of \eqref{eq:morII} and \eqref{eq:reformulation} is an isomorphism in $\SH(\bbV(\mathscr E^0)_\sigma)$.
\end{coro}
We conclude this section by explaining how the main results can be refined to take into account the monodromy action. This is a direct application of the main theorems of \cite{AJM}. 
Let $f:\bbV(\mathscr E)\ra\bbA^1$ be a $\bbG_m$-invariant morphism as before. Then $\bbV(\mathscr E)\times\bbG_m$ is a vector bundle on $X\times\bbG_m$ with an induced linear action of $\bbG_m$ (which is trivial on the factor $\bbG_m$) and since the morphism (see \sectionref{sec:recollection} for the notation)
\[f^{\bbG_m}:\bbV(\mathscr E)\times\bbG_m\ra\bbA^1\]
is $\bbG_m$-invariant, we may apply \theoremref{theo:maintheo} to the morphism $f^{\bbG_m}$ to obtain the following corollary.

\begin{coro}
Let $A\in\SH(\bbV(\mathscr E)_\eta)$ be a $\bbG_m$-equivariant object.
The morphism
\[\acc{\pi}^+_{\sigma!}\acc{e}^{+*}_\sigma\Psi^{\mathrm{mon}}_f(A)\ra\Psi^{\mathrm{mon}}_{f^0}(\pi^+_{\eta!}e^{+*}_\eta A)\]
is an isomorphism in $\QUSH(\bbV(\mathscr E^0)_\sigma)$.
\end{coro}
Similarly by applying \corollaryref{coro:maincoro} to the morphism $f^{\bbG_m}$ we get the following result.
\begin{coro}\label{coro:IImonodromic}
There is an isomorphism in $\QUSH(\bbV(\mathscr E^0)_\sigma)$
\[\acc{\pi}^+_{\sigma!}\acc{e}^{+*}_\sigma\Psi^{\mathrm{mon}}_f(\one_{\bbV(\mathscr E)_\eta})\xra{\simeq}\Thom^{-1}(s^+,\pi^+)|_{\bbV(\mathscr E^0)_\sigma\times\bbG_m}\otimes\Psi^{\mathrm{mon}}_{f^0}(\one_{\bbV(\mathscr E^0)_\eta})\]
\end{coro}
In the above formulation we have used the canonical isomorphism between the inverse Thom motive $\Thom^{-1}(\acc{s}^+,\acc{\pi}^+)$ and the pullback of $\Thom^{-1}(s^+,\pi^+)$ along the projection $\bbV(\mathscr E)\times\bbG_m\ra\bbV(\mathscr E)$ (see \cite[Remarque 1.5.10]{AyoubI}).

\section{Application to virtual motives}\label{sec:virtual}
In this section we relate more precisely \theoremref{theo:maintheo} to the integral identity originally conjectured by Kontsevich and Soibelman in \cite[Section 4.4]{KSDT} (see also \cite{MR2851153}). As in \cite{MR3983293}, we assume that $k$ contains all roots of unity. Let us first recall the stronger form of the integral identity proved by Nicaise and Payne in \cite[Theorem 4.2.1]{MR3983293} as an application of their motivic Fubini theorem for the tropicalization map based on Hrushovski-Kazhdan's theory of motivic volumes of semi-algebraic sets.

\begin{theo}[Nicaise and Payne]\label{theo:NP}
Let $X$ be a $k$-variety. Let $d_1,d_2$ be nonnegative integers and let $\bbG_m$ acts diagonally on $\bbA^{d_1}\times\bbA^{d_2}\times X$
with positive weights on the first factor, negative weights on the second factor and trivially on $X$. Let 
\[f:\bbA^{d_1}\times\bbA^{d_2}\times X\ra\bbA^1\]
be a $\bbG_m$-invariant morphism, $f|_X$ be the restriction of $f$ to $X$. If $x\in X(k)$ is contained in the zero locus of $f|_X$, then the equality
\begin{equation}\label{eq:IIoriginal}
\int_{\bbA^{d_1}}(\psi_f)|_{\bbA^{d_1}}=\bbL^{d_1}\psi_{f|_X,x}
\end{equation}
holds in the ring $\sM^{\hat{\mu}}_k$.
\end{theo}

In \cite[Section 4.4]{KSDT}, equality \eqref{eq:IIoriginal} is conjectured under the assumption that all weights are either 1 or -1 and that the variety $X$ is an affine space $\bbA^{d_3}$ for some positive integer $d_3$ with $x$ taken to be the origin. Using Hrushovski-Kazhdan motivic integration L\^e Quy Thuong proves in \cite[Theorem 1.2]{MR3299104} that this particular case holds in a localization of $\sM^{\hat{\mu}}_k$ obtained by inverting all the elements $1-\bbL^i$ with $i\geq 1$. Note in \cite{MR3983293} this localization is not needed and Nicaise-Payne even show that the equality holds in $\K_0^{\hat{\mu}}(\Var_k)$.

Given \theoremref{theo:maintheo} it is reasonable to expect the following relative form of the integral identity to be true (here we use the notation in  \sectionref{sec:mainresults}).

\begin{conj}\label{conj:IIrelative}
Let $f:\bbV(\mathscr E)\ra\bbA^1$ be a $\bbG_m$-invariant morphism. Let $r$ be the rank of the vector bundle $\bbV(\mathscr E^+)$ over $\bbV(\mathscr E^0)$. Then, the equality 
\begin{equation}\label{eq:IIrelative}
\int_{\pi^+}(\psi_f)|_{\bbV(\mathscr E^+)_\sigma}=\bbL^{r}\psi_{f^0}
\end{equation}
holds in the ring  $\mathscr M^{\hat{\mu}}_{\bbV(\mathscr E^0)_\sigma}$.
\end{conj}

\begin{rema}
Note that \conjectureref{conj:IIrelative} implies \theoremref{theo:NP}. Indeed with the notation of \theoremref{theo:NP}, equality \eqref{eq:IIoriginal} is simply the image of equality \eqref{eq:IIrelative} by the morphim of rings
\[x^*:\mathscr M^{\hat{\mu}}_k\ra \mathscr M^{\hat{\mu}}_{\bbV(\mathscr E^0)_\sigma}.\]
\end{rema}

We now show that \conjectureref{conj:IIrelative} holds if the ring $\mathscr M^{\hat{\mu}}_{\bbV(\mathscr E^0)_\sigma}$ is replaced by the Grothendieck ring of the triangulated category of compact (or constructible) quasi-unipotent motives. This is a direct consequence of \corollaryref{coro:IImonodromic}.
\begin{coro}\label{coro:IIquasiunipotent}
Equality \eqref{eq:IIrelative}:
\[
\int_{\pi^+}(\psi_f)|_{\bbV(\mathscr E^+)_\sigma}=\bbL^{r}\psi_{f^0}
\]
holds in the ring $\K_0(\QUSH_{\ct}(\bbV(\mathscr E^0)_\sigma))$.
\end{coro}

\begin{proof}
By \corollaryref{coro:IImonodromic}, we have an equality in the ring $\K_0(\QUSH_\ct(\bbV(\mathscr E^0)_\sigma))$
\[
\acc{\pi}^+_{\sigma!}\acc{e}^{+*}_\sigma\left[\Psi^{\mathrm{mon}}_f(\one_{\bbV(\mathscr E)_\eta})\right]=\left[ \Thom^{-1}(s^+,\pi^+)|_{\bbV(\mathscr E^0)_\sigma\times\bbG_m}\right]\cdot\left[\Psi^{\mathrm{mon}}_{f^0}(\one_{\bbV(\mathscr E^0)_\eta})\right].\]
Let us observe that the class in the Grothendieck ring of $\Thom^{-1}(s^+,\pi^+)|_{\bbV(\mathscr E^0)_\sigma\times\bbG_m}$ is equal to $[\one_{\bbV(\mathscr E^0)_\sigma\times\bbG_m}(-r)] $ (see \cite[p.1376]{AIS} for details) and thus equal to image of $\bbL^r$ by the ring morphism $\chi_{\bbV(\mathscr E^0)_\sigma,c}$. Using \lemmaref{lemm:chic} and \cite[Theorem 5.3.1]{AJM}, we get
\begin{align*}
\chi_{\bbV(\mathscr E^0)_\sigma,c}\left(\int_{\pi^+}(\psi_f)|_{\bbV(\mathscr E^+)_\sigma}\right)&=\acc{\pi}^+_{\sigma!}\acc{e}^{+*}_\sigma\chi_{\bbV(\mathscr E)_\sigma,c}(\psi_f)
\\
&=\acc{\pi}^+_{\sigma!}\acc{e}^{+*}_\sigma\left[\Psi^{\mathrm{mon}}_f(\one_{\bbV(\mathscr E)_\eta})\right]\\
&=\chi_{\bbV(\mathscr E^0)_\sigma,c}(\bbL^r)\cdot[\Psi^{\mathrm{mon}}_{f^0}(\one_{\bbV(\mathscr E^0)_\eta})]\\
&=\chi_{\bbV(\mathscr E^0)_\sigma,c}(\bbL^r)\cdot \chi_{\bbV(\mathscr E^0)_\sigma,c}(\psi_{f^0})\\
&=\chi_{\bbV(\mathscr E^0)_\sigma,c}(\bbL^r\psi_{f^0}).
\end{align*}
This shows the desired equality.
\end{proof}

\appendix

\renewcommand{\theequation}{\Alph{section}.\arabic{equation}}
\section{Braden's contraction lemma}\label{app:A}
For the sake of completeness, in \appendixref{app:A} and \appendixref{app:B} we give a proof of Braden's hyperbolic localization theorem in motivic homotopy theory in the special case of a vector bundle with a linear action of $\bbG_m$. The reader familiar with \cite{MR1996415} (see also \cite{MR3912059}) will recognize that the proof is exactly the same as Braden's proof originally written in the context of \'etale or transcendental sheaves. Indeed Braden's argument work unchanged in any stable homotopy 2-functor in the sense of \cite[D\'efinition 1.4.1]{AyoubI}. We refer also to \cite{MR3200429} for a different proof of Braden's theorem in the context of constructible sheaves/$\mathscr D$-modules as well as references for applications of Braden's theorem to geometric representation theory.

\begin{defi}\label{defi:equivariant}
Let $S$ be an algebraic $k$-variety endowed with an action of $\bbG_m$ and denote by 
\[a:\bbG_m\times S\ra S\qquad q_S:\bbG_m\times S\ra S\]
the action and the projection. An object $A\in\SH(S)$ is said to be equivariant if $a^*A$ and $q_S^*A$ are isomorphic in $\SH(\bbG_m\times S)$.
\end{defi}
The following remark is a direct consequence of the proper base change theorem and the smooth base change theorem.
\begin{rema}
Let $f:T\ra S$ be a $\bbG_m$-equivariant morphism of $k$-varieties endowed with a $\bbG_m$-action. Then $\bbG_m$-equivariant objects are stable under the four operations $f_*,f^*,f_!,f^!$.
\end{rema}

Let $S$ be an algebraic $k$-variety endowed with an action of $\bbG_m$. The inverse action is the morphism
\[a^{\mathrm{inv}}:\bbG_m\times S\xra{\mathrm{inv}\times \Id}\bbG_m\times S\xra{a} S\]
where $\mathrm{inv}:\bbG_m\ra\bbG_m$ is the inverse morphism of the group scheme $\bbG_m$ and $a:\bbG_m\times S\ra S$ is the given action on $S$. Note that if $A\in\SH(S)$ is weakly equivariant for the given action $a$ on $S$, then it is also equivariant for the inverse action on $S$.

Let $S$ be a smooth $k$-variety and $\mathscr V$ be a locally free $\Osheaf_S$-module. We assume that $\bbG_m$ acts linearly on $\bbV(\mathscr V)$ and consider the associated decomposition $\mathscr V=\oplus_{n\in \bbZ}\mathscr V_n$. Given an integer $r\in\bbZ $, we set $\mathscr V^{<r}=\oplus_{n<r}\mathscr V_n$ and $\mathscr V^{\geq r}=\oplus_{n\geq r}\mathscr V_n$. Consider the natural morphisms
\begin{equation}\label{eq:diagramcontraction}
\xymatrix{{\bbP(\mathscr V^{<r})}\ar[r]^-{e_r}\ar@/^2em/[rr]^{\Id_{\bbP(\mathscr V^{<r})}}& {\bbP(\mathscr V)\setminus\bbP(\mathscr V^{\geq r})}\ar[r]^-{\pi_r} &{\bbP(\mathscr V^{<r})}}
\end{equation}
and the induced natural transformation
\begin{equation}\label{eq:morcontraction}
\pi_{r*}\ra \pi_{r*}e_{r*}e_r^*\simeq e_r^*.
\end{equation}
Let $\tau_r:\bbP(\mathscr V^{<r})\ra S$ be the structural morphism. 
We have Braden's contraction lemma (see \cite[Lemma 6]{MR1996415}).
\begin{prop}\label{prop:cont}
Let $r\in\bbZ$ be an integer and $A\in\SH(\bbP(\mathscr V)\setminus \bbP(\mathscr V^{\geq r}))$ be a $\bbG_m$-equivariant object . Then, the morphism 
\[\tau_{r*}\pi_{r*}A\ra\tau_{r*}e_r^*A\]
is an isomorphism in $\SH(S)$.  
\end{prop}

\begin{proof}
As $r$ is fixed in the proof, to shorten the notation we set $e:=e_r$, $\pi:=\pi_r$, $\tau:=\tau_r$ and  $T:=\bbP(\mathscr V)\setminus\bbP(\mathscr V^{\geq r})$. We denote by $u:U\hookrightarrow T$ the open immersion of the complement of $\bbP(\mathscr V^{<r})$ in $T$. The given linear action of $\bbG_m$ on $\bbV(\mathscr V)$ induces an action of $\bbG_m$ on $T$. Let $\Gamma$ be the schematic closure in $\bbA^1\times T\times T$ of the graph of the action morphism $a:\bbG_m\times T\ra T$.
Denote by $p_1,p_2:\Gamma\ra\bbA^1\times T$ the morphisms obtained by composing the closed immersion of $\Gamma$ into with $\bbA^1\times T\times T$ respectively with the projection onto the first two factors and with the projection onto the first and third factors. We recall the following facts:
\begin{enumerate}
\item[(a)] the morphism obtained from $p_2$ by base change along the open immersion $\bbG_m\times T\hookrightarrow \bbA^1\times T$ is an isomorphism;
\item[(b)] the morphism $p_1$ is proper;
\item[(c)] if $u':\bbG_m\times U\ra \Gamma$ is the morphism obtained from the graph morphism of the morphism $\bbG_m\times U\ra U$ induced by the action, then the commutative square
\begin{equation}\label{eq:squarecont}
\xymatrix{{\bbG_m\times U}\ar[d]^-{a}\ar[r]^-{u'} & {\Gamma}\ar[d]^-{q_2\circ p_2}\\
{U}\ar[r]^-{u} & {T}}
\end{equation}
is cartesian.
\end{enumerate}
A proof of these properties can be found in \cite[Lemma 2.14, p.272]{MR3912059} (see also \cite[Proof of Lemma 6]{MR1996415}).

As a consequence of the localization triangle $u_!u^*A\ra A\ra e_*e^*A$, to prove the proposition, it suffices to show the vanishing of the motive $M:=\tau_*\pi_*u_!u^*A$. By applying the smooth base change theorem to the cartesian square
\[\xymatrix{
{\bbA^1\times T}\ar[r]^-{q_2}\ar[d]^-{q_1}\ar@{}[rd]|{\square} &{T}\ar[d]^-{\tau\circ\pi}\\
{\bbA^1\times S}\ar[r]^-{p_S} & {S}}\]
we get that $p_S^*M$ is isomorphic to $q_{1*}q_2^*u_!u^*A$ and we may consider the morphism in $\SH(\bbA^1\times S)$
\begin{equation}\label{eq:cont1}
p_S^*M\ra N:=q_{1*}p_{2*}p_2^*q_2^*u_!u^*A
\end{equation}
induced by the natural transformation $\Id\ra p_{2*}p_2^*$. The morphism \eqref{eq:cont1} has the following properties:
\begin{enumerate}
\item the morphism $j_S^*p_S^*M\ra j_S^*N$ is an isomorphism in $\SH(\bbG_m\times S)$;
\item the composition 
\begin{equation}\label{eq:cont2}
\xymatrix{{j_{S!}j_S^*p^*S}\ar[r] & {p_S^*M}\ar[r]^-{\eqref{eq:cont1}} & {N}}
\end{equation}
is an isomorphism in $\SH(\bbA^1\times S)$.
\end{enumerate}
Note that (1) follows directly from (a). Let us prove (2). The property (1) ensures that the lower horizontal morphism in the commutative square 
\[\xymatrix{{p_S^*M}\ar[r]^-{\eqref{eq:cont1}} & {N}\\
{j_{S!}j_S^*p_S^*M}\ar[r]\ar[u] & {j_{S!}j_S^*N}\ar[u]}\]
is an isomorphism. Hence, to prove (2), it suffices to show that the morphism $j_{S!}j_S^*N\ra N$ is an isomorphism. Therefore, we have to show the vanishing of $i_S^*N$.

By (b) and the fact that $q_1\circ p_2=q_1\circ p_1$, the motive $N$ is isomorphic to $q_{1*}p_{1!}p_2^*q_2^*u_!u^*A$. Using the proper base change theorem in the square \eqref{eq:squarecont}, which is cartesian by (c), we get an isomorphism between $N$ and $q_{1*}p_{1!}u'_!a^*u^*A$. Using the fact that $A$ is assumed to be $\bbG_m$-equivariant and observing that $p_1\circ u'=j\times u$ we get an isomorphism between $N$ and $q_{1*}(j\times u)_!q_U^*u^*A$. Using the commutative diagram
\[\xymatrix{{\bbG_m\times U}\ar[r]^-{\Id_{\bbG_m}\times u}\ar[d]^-{q_U}\ar@{}[rd]|{\square}\ar@/^2em/[rr]^{j\times u} & {\bbG_m\times T}\ar[r]^-{j_T}\ar[d]^-{q_T} & {\bbA^1\times T}\ar[ld]^-{q_2}\\
{U}\ar[r]^-{u} & {T} & {}}\]
and the proper base change theorem, we finally see that $N$ is isomorphic to the motive $q_{1*}j_{T!}j_T^*q_2^*u_!u^*A$. To prove the statement, it suffices to prove that the morphism
\begin{equation}\label{eq:cont4}
q_{1*}q_2^*u_!u^*A\ra q_{1*}i_{T*}i_T^*q_2^*u_!u^*A
\end{equation}
becomes an isomorphism after applying the functor $i^*_S$. Note that $q_1\circ i_T=i_S\circ (\tau\circ\pi)$ and $q_2\circ i_T=\Id_T$ so that the right-hand side of \eqref{eq:cont4} is isomorphic $i_{S*}M$. On the other hand, by the smooth base change, the left-hand side of \eqref{eq:cont4} is isomorphic to $p_S^*M$. With these identifications, the morphism \eqref{eq:cont4} is given by the canonical morphism
\[p_S^*M\ra i_{S*}i_S^*p_S^*M=i_{S*}M.\]
Since this morphism becomes an isomorphism after applying the functor $i_S^*$ the property (2) is proved.

To finish the proof of \propositionref{prop:cont}, consider now the endomorphism of $p_S^*M$ given by the composition
\begin{equation}\label{eq:cont3}
\xymatrix{{p_S^*M}\ar[r]^-{\eqref{eq:cont1}} & {N}\ar[r]^-{\eqref{eq:cont2}^{-1}} & {j_{S!}j_S^*p_S^*M}\ar[r] & {p_S^*M.}}
\end{equation}
Let us denote by $1_S:S\hookrightarrow \bbG_m\times S$ the unit section. Then,
we have commutative square
\[\xymatrix{{p_{S*}1_{S*}1_S^*p_S^*M}\ar[d]^-{\simeq} & {p_{S*}p_S^*M}\ar[l]\ar[r]\ar[d]^-{p_{S*}\eqref{eq:cont3}} & {p_{S*}i_{S*}i_S^*p_S^*M}\ar[d]^-{0}\\
 {p_{S*}1_{S*}1_S^*p_S^*M} & {p_{S*}p_S^*M}\ar[r]\ar[l] & {p_{S*}i_{S*}i_S^*p_S^*M}}\]
This shows that $p_{S*}p_S^*M$ vanishes and therefore that $M$ vanishes since by homotopy invariance the canonical morphism $M\ra p_{S*}p_S^*M$ is an isomorphism.
\end{proof}

\section{Braden's hyperbolic localization theorem}\label{app:B}

We use the notation introduced in \sectionref{sec:mainresults} in particular the notation from diagram \eqref{eq:diagnot}. Since $\pi^-\circ s^-=\Id_{\bbV(\mathscr E^0)}$ and $\pi^+\circ s^+=\Id_{\bbV(\mathscr E^0)} $, we have natural transformations 
\[\pi^-_*\ra \pi^-_*s^-_*s^{-*}\xra{\simeq}s^{-*},\qquad s^{+!}\xra{\simeq} \pi^+_!s^+_!s^{+!}\ra \pi^+_!.\]
As a consequence of \propositionref{prop:cont}, we have the following result.
\begin{prop}\label{prop:simplecase}
Let $A\in\SH(\bbV(\mathscr E^-))$ and $B\in\SH(\bbV(\mathscr E^+))$ be $\bbG_m$-equivariant objects. Then, the natural morphisms
\[\pi^-_*A\ra s^{-*}A\qquad\textrm{and}\qquad s^{+!}B\ra \pi^+_!B\]
are isomorphisms in $\SH(\bbV(\mathscr E^0))$.
\end{prop}

\begin{proof}
We only have to prove that the morphism $\pi^-_*A\ra s^{-*}A$ is an isomorphism since duality then implies that the other one is an isomorphism too. As we may replace $X$ by $\bbV(\mathscr E^0)$ as base scheme, we may assume that $\mathscr E^0=0$. We consider on $\bbV(\mathscr E)$ the inverse action of $\bbG_m$ and consider the $\Osheaf_X$-module $\mathscr V=\Osheaf_X\oplus \mathscr E^-$ where $\Osheaf_X$ has weight zero. The induced decomposition on $\mathscr V$ is given by 
\[\mathscr V_n=\begin{cases} \mathscr E_{-n} & \textrm{if $n>0$}\\\Osheaf_X & \textrm{if $n=0$}\\0 & \textrm{if $n<0 $} \end{cases}\]
so that $\mathscr V^{<1}=\Osheaf_X$ and $\mathscr V^{\geq 1}=\mathscr E^-$. Thus we have $\bbP(\mathscr V^{<1})=X$ as well as $\bbP(\mathscr V)\setminus\bbP(\mathscr V^{\geq 1})=\bbV(\mathscr E^-)$, and the result follows from \propositionref{prop:cont}.
\end{proof}

Consider now, given an object $A$ in $\SH(\bbV(\mathscr E))$, the morphism
\begin{equation}\label{eq:morpurity}
s^{-*}e^{-!}A\ra s^{+!}e^{+*}A
\end{equation}
obtained as the composition
\[s^{-*}e^{-!}A\ra s^{-*}e^{-!}e^+_*e^{+*}A\ra s^{-*}s^-_*s^{+!}e^{+*}A\ra s^{+!}e^{+*}A\]
where the first morphism is induced by the unit of the adjunction $e^{+*}\dashv e^+_*$, the second morphism is induced by the inverse of the exchange 2-isomorphism $s^-_*s^{+!}\ra e^{-!}e^+_*$, and the third one is induced by the counit of the adjunction $s^{-*}\dashv s^-_*$.
\begin{rema}\label{rema:twoiso}
Note that in the composition that defines \eqref{eq:morpurity} the second and the third morphisms are isomorphisms (for the last one it follows from the fact that $s^-$ is a closed immersion).
\end{rema}

\begin{prop}\label{prop:purity}
Let $A\in\SH(\bbV(\mathscr E))$ be a $\bbG_m$-equivariant object. Then, the morphism 
\[s^{-*}e^{-!}A\ra s^{+!}e^{+*}A\]
defined in \eqref{eq:morpurity} is an isomorphism in $\SH(\bbV(\mathscr E^0))$.
\end{prop}

Before proving \propositionref{prop:purity}, let us mention that this result can be interpreted in terms of the generalization of absolute purity introduced by Ayoub in \cite{AyoubEtale}. The morphism \eqref{eq:morpurity} as defined above is the one considered in \cite{MR1996415,MR3200429,MR3912059} in their respective context. However it coincides also with the exchange morphism associated with the cartesian squares of closed immersions in \eqref{eq:diagnot} as considered by Ayoub in \cite[D\'efinition 7.1]{AyoubEtale} (this can be checked easily from the definitions of the various exchange morphisms given in \cite[Proposition 1.4.15 and page 80]{AyoubI}). Hence, \propositionref{prop:purity} means that $\bbG_m$-equivariant objects in $\SH(\bbV(\mathscr E))$ are pure relatively to the cartesian square in \eqref{eq:diagnot} in the sense introduced by Ayoub in \cite[D\'efinition 7.1]{AyoubEtale}.

\begin{proof}[Proof of \propositionref{prop:purity}]
As we may replace $X$ by $\bbV(\mathscr E^0)$ as base scheme, we may assume that $\mathscr E^0=0$. If either $\mathscr E^+$ or $\mathscr E^-$ is zero the proposition holds true. Hence we may also assume  that $\mathscr E^+$ and $\mathscr E^-$ are non zero. Consider the closed immersion $e^+:\bbV(\mathscr E^+)\hookrightarrow \bbV(\mathscr E)$, the open complement $U=\bbV(\mathscr E)\setminus \bbV(\mathscr E^+)$  and denote by $u:U\hookrightarrow\bbV(\mathscr E)$ the open immersion. Consider the localization triangle
\[u_!u^*A\ra A\ra e^+_*e^{+*}A\xra{+1}.\]
By construction of the morphism \eqref{eq:morpurity} and \remarkref{rema:twoiso}, we have to prove that $s^{-*}e^{-!}u_!u^*A$ vanishes. By \propositionref{prop:simplecase}, the motive $s^{-*}e^{-!}u_!u^*A$ is isomorphic to $\pi^-_*e^{-!}u_!u^*A$ and thus we have to prove that $\pi^-_*e^{-!}u_!u^*A$ vanishes. We consider on $\bbV(\mathscr E)$ the inverse action of $\bbG_m$ and consider the $\Osheaf_X$-module $\mathscr V=\Osheaf_X\oplus \mathscr E$ where $\Osheaf_X$ has weight zero. The induced decomposition on $\mathscr V$ is given by 
\[\mathscr V_n=\begin{cases} \mathscr E_{-n} & \textrm{if $n\neq 0$}\\\Osheaf_X & \textrm{if $n=0$}\end{cases}\]
so that $\mathscr V^{<1}=\Osheaf_X\oplus\mathscr E^+$, $\mathscr V^{\geq 1}=\mathscr E^-$ and $\mathscr V^{<0}=\mathscr E^+$, $\mathscr V^{\geq 0}=\Osheaf_X\oplus \mathscr E^-$. Then, we have the commutative diagram
\[\xymatrix{{\bbV(\mathscr E^+)}\ar[d]^-{\rho^+}\ar[r]^-{e^+} & {\bbV(\mathscr E)}\ar[d]^-{\rho} & {\bbV(\mathscr E^-)}\ar[l]_-{e^-}\ar[rd]^-{\pi^-}\\
{\bbP(\mathscr V^{<1})}\ar[r]^-{e_1} & {\bbP(\mathscr V)\setminus\bbP(\mathscr V^{\geq 1})}\ar[r]^-{\pi_1} & {\bbP(\mathscr V^{<1})}\ar[r]^-{\tau_1} & {X}}\]
in which $\rho$ is the usual open immersion and $\rho^+$ the induced morphism. Note that $\rho\circ e^-$ is a closed immersion. Let $V$ be the open complement of $\bbV(\mathscr E^-)$ in $\bbP(\mathscr V)\setminus\bbP(\mathscr V^{\geq 1})$. Note that $V=\bbP(\mathscr V)\setminus\bbP(\mathscr V^{\geq 0})$. Denote by $v:V\hookrightarrow \bbP(\mathscr V)\setminus\bbP(\mathscr V^{\geq 1})$ the open immersion.
Let $B=\rho_!u_!u^*A$ so that $\rho^!B=u_!u^*A$. Consider the localization exact triangle
\[\rho_*e^-_*e^{-!}\rho^!B\ra B\ra v_*v^*B\xra{+1}\]
and apply to it the functor $\tau_{1*}\pi_{1*}$ to obtain an exact triangle
\[\pi^-_*e^{-!}u_!u^*A\ra \tau_{1*}\pi_{1*}B\ra \tau_{1*}\pi_{1*}v_*v^*B\xra{+1}.\]
To prove the proposition, it suffices now to show that both $\tau_{1*}\pi_{1*}B$ and $\tau_{1*}\pi_{1*}v_*v^*B$ vanish. Let us first show that $\tau_{1*}\pi_{1*}B=0$ by considering the commutative diagram
\[\xymatrix{{\bbV(\mathscr E^+)}\ar[r]^-{e^+}\ar[d]^-{\rho^+}\ar@{}[rd]|{\square} & {\bbV(\mathscr E)}\ar[d]^-{\rho}\ar[d] &{U=\bbV(\mathscr E)\setminus \bbV(\mathscr E^+)}\ar[d]\ar[l]_-u\\
{\bbP(\mathscr V^{<1})}\ar[r]^-{e_1} & {\bbP(\mathscr V)\setminus\bbP(\mathscr V^{\geq 1})} & {\bbP(\mathscr V)\setminus \bbP(\mathscr V^{<1}).}\ar[l]}\]
Since $B$ is $\bbG_m$-equivariant, we know by \propositionref{prop:cont} that $\tau_{1*}\pi_{1*}B$ is isomorphic to $\tau_{1*}e^*_1B$ and the results follows from the fact that $e_1^*B=\rho^+_!e^{+*}u_!A=0$. It remains to prove that $\tau_{1*}\pi_{1*}v_*v^*B=0$. Consider the commutative diagram
\begin{equation}\label{eq:diapurityend}
\xymatrix@R=.5cm{{\bbV(\mathscr E)\setminus \bbV(\mathscr E^+)}\ar[d]^-{u}\ar@{}[rrdd]|{\square} & &{\bbV(\mathscr E)\setminus(\bbV(\mathscr E^+)\cup\bbV(\mathscr E^-))}\ar[ll]_-{w}\ar[dd]^-{u'}& {\emptyset}\ar[dd]\ar[l]\ar@{}[ldd]|{\square}\\
 {\bbV(\mathscr E)}\ar[d]^-{\rho}& {}\\
 {\bbP(\mathscr V)\setminus \bbP(\mathscr V^{\geq 1})}\ar[d]^-{\pi_1} & & {V=\bbP(\mathscr V)\setminus\bbP(\mathscr V^{\geq 0})}\ar[ll]_-{v}\ar[d]^-{\pi_0} & {\bbP(\mathscr V^{<0})}\ar[l]_-{e_0}\\
 {\bbP(\mathscr V^{<1})}\ar[r]^-{\tau_1}& {X} & {\bbP(\mathscr V^{<0}).}\ar[l]_-{\tau_0}}
 \end{equation}
Applying the proper base change theorem to the cartesian square on the left in \eqref{eq:diapurityend}, we see that $\tau_{1*}\pi_{1*}v_*v^*B$ is isomorphic to $\tau_{0*}\pi_{0*}u'_!w^*u^*A$ and therefore to $\tau_{0*}e_0^*u'_!w^*u^*A $ by \propositionref{prop:cont}. Now it suffices to apply the proper base change theorem to the cartesian square on the right-hand side of \eqref{eq:diapurityend} to get the desired vanishing.
\end{proof}

Braden's hyperbolic localization theorem in motivic homotopy theory is then a direct consequence of \propositionref{prop:purity} and \propositionref{prop:simplecase}:

\begin{theo}\label{theo:BradenIso}
Let $A\in\SH(\bbV(\mathscr E))$ be a $\bbG_m$-equivariant object. Then, the morphism
\begin{equation}\label{eq:morHypLoc}
\pi^-_*e^{-!}A\ra \pi^+_!e^{+*}A
\end{equation}
is an isomorphism in $\SH(\bbV(\mathscr E^0))$.
\end{theo}
Note that the morphism \eqref{eq:morHypLoc} is nothing but the morphism \eqref{eq:morpurity} composed with the isomorphisms of \propositionref{prop:simplecase}.

\bibliographystyle{amsplain}
\bibliography{HypLoc}

\end{document}